\documentclass[11 pt,english]{article}
\usepackage[activeacute]{babel}
\usepackage{amsfonts} 
\usepackage{mathrsfs} 
\usepackage{amssymb, amsmath, amsfonts}
\usepackage{makeidx}
\usepackage{epsfig}
\usepackage{color}
\usepackage[T1]{fontenc}
\usepackage{graphics,graphicx,graphpap}
\usepackage[ansinew]{inputenc}
\usepackage{amsthm}
\usepackage{subfigure}
\usepackage{tikz}

\usetikzlibrary{calc, positioning, through, math, arrows, backgrounds}
\tikzstyle{vtx}=[inner sep=1pt,draw, shape=circle, font=\tiny]
\tikzstyle{line}=[inner sep=3pt,draw, shape=rectangle, line width = 3pt]
\tikzset{>=stealth}
\tikzstyle{lbl}=[inner sep = 1 pt, fill = white, midway]

\newcommand\blfootnote[1]{%
  \begingroup
  \renewcommand\thefootnote{}\footnote{#1}%
  \addtocounter{footnote}{-1}%
  \endgroup
}

\if@twoside \oddsidemargin 6pt \evensidemargin 6pt \marginparwidth 20pt
 \else \oddsidemargin 6pt \evensidemargin 6pt \marginparwidth 20pt \fi
\newcommand{\ZZ}{\mathbb{Z}}

\newcommand{\G}{\Gamma}
\newcommand{\Cay}{\mathop{\rm Cay}}

\newcommand{\la}{\langle}
\newcommand{\ra}{\rangle}

\newcommand{\Circ}{\mathop{\rm Circ}}

\newtheorem{theorem}{Theorem}[section]
\newtheorem{proposition}[theorem]{Proposition}
\newtheorem{corollary}[theorem]{Corollary}
\newtheorem{lemma}[theorem]{Lemma}
\newtheorem{question}[theorem]{Question}

\theoremstyle{definition}

\newtheorem{assumption}[theorem]{Assumption}

\newcommand{\Claim}{{Claim~}}

\begin{document}

\begin{center}
\Large{\textbf{On the fractional matching extendability of Cayley graphs of Abelian groups}} \\ [+4ex]
\end{center}

\begin{center}
Bo\v stjan Kuzman{\small $^{a}$},\   
Primo\v z \v Sparl{\small $^{a,b,c,*}$}
\\

\medskip
{\it {\small
$^a$ University of Ljubljana, Faculty of Education, Ljubljana, Slovenia\\
$^b$University of Primorska, Institute Andrej Maru\v si\v c, Koper, Slovenia\\
$^c$Institute of Mathematics, Physics and Mechanics, Ljubljana, Slovenia\\
}}
\end{center}

\blfootnote{
Email addresses: 
bostjan.kuzman@pef.uni-lj.si (B.~Kuzman), 
primoz.sparl@pef.uni-lj.si (P.~\v Sparl)\\
* - corresponding author
}


\hrule

\begin{abstract}
Fractional matching extendability is a concept that brings together two widely studied topics in graph theory, namely that of fractional matchings and that of matching extendability. A {\em fractional matching} of a graph $\Gamma$ with edge set $E$ is a function $f$ from $E$ to the real interval $[0,1]$ with the property that for each vertex $v$ of $\Gamma$, the sum of $f$-values of all the edges incident to $v$ is at most $1$. When this sum equals $1$ for each vertex $v$, the fractional matching is {\em perfect}. A graph of order at least $2t+1$ is {\em fractional $t$-extendable} if it contains a matching of size $t$ and if each such matching $M$ can be extended to a fractional perfect matching in the sense that the corresponding function $f$ assigns value $1$ to each edge of $M$. 

In this paper, we study fractional matching extendability of Cayley graphs of Abelian groups. We show that, except for the odd cycles, all connected Cayley graphs of Abelian groups are fractional $1$-extendable and we classify the fractional $2$-extendable Cayley graphs of Abelian groups. This extends the classification of $2$-extendable (in the classical sense) connected Cayley graphs of Abelian groups of even order from 1995, obtained by Chan, Chen and Yu.
\end{abstract}
\hrule

\begin{quotation}
\end{quotation}

\section{Introduction}
\label{sec:Intro}

In~\cite{SchUllbook}, Scheinerman and Ullman presented various fascinating aspects of fractional graph theory and showed that by allowing certain parameters to be real numbers and not only integers, one obtains interesting generalizations of classical concepts such as the independence number, minimum vertex cover, packing number, etc. The corresponding results generalize their classical counterparts and can thereby give bounds on certain well-studied parameters of graphs. The fractional parameters are often much easier to compute. This way, efficient approximation algorithms for some classical 
NP-hard or NP-complete problems can be obtained.

One example, where the fractional generalization is in particular fruitful, is that of matchings in graphs. Instead of choosing a set of independent edges, in the fractional setting one assigns certain values from the real interval $[0,1]$ to the edges of the graph in such a way that the sum of the values at each vertex is at most $1$ (see Section~\ref{sec:Prelim} for a precise definition and for the definitions of all other terms not defined in the Introduction). When this sum equals $1$ for all the vertices, one obtains what is called a fractional perfect matching. Fractional matchings in (hyper)graphs have been studied extensively for more than four decades (see~\cite{Bal81, MaQiaShi22, O16, Pul87, YanYua22} for a few examples). In particular, the problem of existence of fractional perfect matchings in (hyper)graphs is a common topic (see~\cite{DevKah17, PanLiu23} for two of the more recent papers on the subject).

A related concept is that of matching extendability, where one is interested in extending a given matching to a perfect matching. Introduced by Plummer in 1980~\cite{Plu80}, the concept has attracted much attention since (see for instance the surveys~\cite{Plu94, Plu08} and~\cite{CioKooLi17, CioLi14, KutMarMikSpa19, MikSpa09, ZhavDa23} for some more recent results). We mention that the corresponding matching extension problem  in  general graphs is also difficult, as proven recently in~\cite{HacKos18} (see also~\cite{Plu08}).

It is thus natural to investigate the concepts at the intersection of these two topics - the fractional (perfect) matchings and extendability. This leads to the so-called fractional extendability of matchings, first considered by Ma~\cite{MaPhD} and later also investigated in~\cite{neki-kinezi, YuCao08}. Here the question is, whether in a given graph, one can arbitrarily choose a prescribed number of independent edges and extend the chosen set to a fractional perfect matching (see Section~\ref{sec:Prelim} for a precise definition). This concept is the main theme of the paper. However, certain restrictions are needed in order to obtain some classification results. In this paper, we focus on a well-studied class of regular graphs, namely, Cayley graphs of Abelian groups. We mention that for this class of graphs,  extendability of matchings in the classical sense has been previously investigated (see~\cite{ChaCheYu95, MikSpa09}).
\medskip

The paper is organized as follows. In Section~\ref{sec:Prelim}, we gather several definitions and known results that we will need in the subsequent sections. In Section~\ref{sec:CayAb}, we study the fractional extendability for Cayley graphs of Abelian groups. In Subsection~\ref{subsec:1-ext}, we classify the fractional $1$-extendable ones (Corollary~\ref{cor:CayAb1fext}). The main result of Section~\ref{sec:CayAb}, however, is the complete classification of all fractional $2$-extendable connected Cayley graphs of Abelian groups (Theorem~\ref{the:CayAb2}). In Subsection~\ref{subsec:even}, the case of graphs of even order is settled using results from~\cite{ChaCheYu95}. There are no analogous results for graphs of odd order, and so the analysis of these is much more complex. The series of intermediate results that altogether lead to our main result, is carried out in Subsection~\ref{ssec: f2e-odd}. Finally, some interesting questions and possible directions for future research are given in Section~\ref{sec:conclude}.

\section{Preliminaries}
\label{sec:Prelim}

Throughout this paper all graphs are assumed to be simple and undirected. Let $\G = (V, E)$ be a graph with vertex set $V$ and edge set $E$. An edge $e \in E$ is an unordered pair of distinct (adjacent) vertices of $\G$. Therefore, for an edge $e \in E$ and a vertex $v \in V$, notation $v \in e$ stands for ``$v$ is an endvertex of $e$''. We also use $N(v)$ for the neighbourhood of a vertex $v \in V$. For a subset $U \subset V$, we denote by $\G[U]$ the subgraph of $\G$ induced on $U$, and by $\G-U$ the subgraph of $\G$ induced on $V \setminus U$ (that is, $\G - U = \G[V \setminus U])$. 

A {\em fractional matching} of a graph $\G = (V, E)$ is a function $f$ from $E$ to the closed interval $[0,1] \subset \mathbb{R}$ such that $\sum_{v \in e}f(e) \leq 1$ for each vertex $v$ of $\G$. In the case that $\sum_{v \in e}f(e) = 1$ holds for each $v \in V$, the function $f$ is said to be a {\em fractional perfect matching} of $\G$. 
Note that each matching of $\G$ is a fractional matching of $\G$ and each perfect matching of $\G$ is a fractional perfect matching of $\G$. However, while an odd cycle (a cycle of odd length) does not have a perfect matching, it does have a fractional perfect matching (one simply assigns $1/2$ to each edge).

It follows from~\cite[Theorem~2.1.5]{SchUllbook}, that if a graph $\G$ admits a fractional perfect matching, then one can in fact find a fractional perfect matching $f$ of $\G$, for which $f(e) \in \{0, 1/2, 1\}$ holds for each edge $e$ of $\G$. Clearly, the set of edges having nonzero value of $f$ then consists of a set of independent edges (with $f$-value $1$) and a set of cycles (whose edges have $f$-value $1/2$). This gives the following useful corollary (see~\cite[Proposition~2.2.2]{SchUllbook}).

\begin{proposition}{\rm \cite{SchUllbook}}
\label{pro:characterization}
Let $\G = (V, E)$ be a graph. Then $\G$ admits a fractional perfect matching if and only if it has a spanning subgraph whose connected components are either edges or odd cycles.
\end{proposition}

The following nice generalization of the celebrated result of Tutte~\cite{Tut47} on the existence of perfect matchings in graphs gives a necessary and sufficient condition for the existence of a fractional perfect matching in a graph. For our purposes, a straightforward corollary of this result will be very useful.

\begin{proposition}{\rm \cite[Theorem~2.2.4]{SchUllbook}}
\label{prop:iso}
Let $\G = (V,E)$ be a graph. Then $\G$ has a fractional perfect matching if and only if for each subset $U \subseteq V$ the number of isolated vertices of the graph $\G - U$ is at most $|U|$.
\end{proposition}

\begin{corollary}
\label{cor:bip}
Let $\G = (V,E)$ be a graph containing an independent set $U \subset V$ such that $|U| > |V \setminus U|$. Then $\G$ has no fractional perfect matching.
\end{corollary}

We now review the concept of fractional extendability of matchings, which was introduced in~\cite{MaPhD} (but see also~\cite{neki-kinezi, YuCao08}). For a graph $\G$ and a non-negative integer $t$, the graph $\G$ is said to be {\em fractional $t$-extendable}, if it is of order at least $2t+1$, admits a matching of size $t$, and if each matching $M$ of $\G$ of size $t$ can be extended to a fractional perfect matching of $\G$ in the sense that the corresponding function $f$ satisfies $f(e) = 1$ for each $e \in M$.

The following corollary of Proposition~\ref{pro:characterization} characterizes fractional $t$-extendable graphs and will be used throughout Section~\ref{sec:CayAb}. 

\begin{proposition}
\label{pro:kriterij}
Let $t \geq 1$ be an integer and let $\G$ be a graph of order at least $2t+1$ which admits a matching of size $t$. Then $\G$ is fractional $t$-extendable if and only if for each matching $M$ of size $t$, the subgraph of $\G$ obtained by removing the endvertices of the edges of $M$, has a spanning subgraph whose connected components are either edges or odd cycles.
\end{proposition} 

Let us point out that while the classical concept of extendability of matchings is only interesting for graphs of even orders, the concept of fractional extendability can be studied for all graphs. Moreover, when considering only graphs of even order, it is of course clear that each $t$-extendable graph (a graph in which each matching of size $t$ can be extended to a perfect matching) is fractional $t$-extendable. But the converse does not hold. For example, let $\G$ be the graph obtained by taking two copies of the complete graph $K_4$ and adding an edge $e$, joining one vertex from one copy to one vertex from the other copy. The matching $M = \{e\}$ clearly cannot be extended to a perfect matching of $\G$. On the other hand, Proposition~\ref{pro:kriterij} implies that it can be extended to a fractional perfect matching of $\G$. It is not difficult to see that all other matchings of size $1$ can also be extended to a fractional perfect matching of $\G$, and so $\G$ is fractional $1$-extendable. This shows that the condition of fractional $t$-extendability is weaker than that of $t$-extendability.



We conclude this section by recalling the well-known notion of a Cayley graph. Let $G$ be a group and let $S \subset G$ be an inverse-closed subset not containing the identity. Then the {\em Cayley graph} $\Cay(G ; S)$ of $G$ with respect to the {\em connection set} $S$ is the graph with vertex set $G$, in which $N(g) = \{gs \colon s \in S\}$ for each $g \in G$. Note that the graph $\Cay(G ; S)$ is regular of valency $|S|$ and is connected if and only if $\la S \ra = G$. In the case that $G$ is the (additive) cyclic group $\ZZ_n = \{0,1, \ldots , n-1\}$ (with addition modulo $n$), we denote $\Cay(\ZZ_n ; S)$ by $\Circ(n ; S)$ and say that the graph is a {\em circulant}.

Throughout the next section we will be working with Cayley graphs of Abelian groups. When working with specific Abelian groups (like cyclic groups or their direct products), we will assume that the corresponding operation is additive (and denote the corresponding groups by $\ZZ_n$, $\ZZ_n \times  \ZZ_m$, etc.). However, when working with a general Abelian group $A$, it will be more convenient to assume that the corresponding operation is multiplicative. This should cause no confusion.




\section{Cayley graphs of Abelian groups}
\label{sec:CayAb}

As mentioned in the Introduction, the classical extendability of matchings for Cayley graphs of Abelian groups was first studied in~\cite{ChaCheYu95}, where the $2$-extendable Cayley graphs of Abelian groups (of even order) were classified. In this section, we consider fractional extendability of connected Cayley graphs of Abelian groups. We first show that, with the exception of odd cycles, all such graphs are fractional $1$-extendable (see Corollary~\ref{cor:CayAb1fext}). We then classify the fractional $2$-extendable examples (see Theorem~\ref{the:CayAb2}). 

\subsection{Fractional $1$-extendability}
\label{subsec:1-ext}

Clearly, odd cycles are not fractional $1$-extendable. In this subsection, we show that all other connected Cayley graphs of Abelian groups of order at least $3$ are fractional $1$-extendable. Our proof is based on the following result from~\cite{MikSpa12} on extendability of paths to Hamilton cycles (cycles containing all vertices of the graph). 

\begin{proposition}{\rm \cite[Theorem~4.1]{MikSpa12}}
\label{pro:3HCext}
Let $\G$ be a connected Cayley graph of an Abelian group of order at least $4$. Then each path of length $3$ of $\G$ lies in a Hamilton cycle of $\G$ if and only if $\G$ is not isomorphic to one of the following graphs:
\begin{itemize}
\itemsep = 0pt
\item[{\rm (i)}] $\Circ(4m ; \{\pm 1, 2m\})$, $m \geq 2$,
\item[{\rm (ii)}] $\Circ(4m+2 ; \{\pm 2, 2m+1\})$, $m \geq 1$,
\item[{\rm (iii)}] $\Circ(4m+2 ; \{\pm 1, \pm 2m\})$, $m \geq 1$,
\item[{\rm (iv)}] $\Circ(2m+1 ; \{\pm 1, \pm 3\})$, $m \geq 3$,
\item[{\rm (v)}] $\Circ(m ; \{\pm 1, \pm 2\})$, $m \geq 6$.
\end{itemize}
\end{proposition}

We mention that the graphs from the first two families of Proposition~\ref{pro:3HCext}, are the well-known M\"obius ladders and prisms, respectively. The ones from family (iii) are often called the wreath graphs.

\begin{corollary}
\label{cor:CayAb1fext}
Let $\G$ be a connected Cayley graph of an Abelian group of order at least $3$. Then $\G$ is fractional $1$-extendable if and only if it is not an odd cycle. 
\end{corollary}

\begin{proof}
As already mentioned, odd cycles are not fractional $1$-extendable. The even cycles clearly are (they are in fact $1$-extendable). We can thus assume that $\G$ is of valency at least $3$.

Denote the corresponding Abelian group (with multiplicative operation) and connection set by $A$ and $S$, respectively, that is, $\G = \Cay(A ; S)$. Suppose first that $\G$ is not isomorphic to a member of one of the five families of circulants from Proposition~\ref{pro:3HCext}. Recall that each Cayley graph is vertex-transitive, meaning that the automorphism group is transitive on the vertex set. To show that $\G$ is fractional $1$-extendable, it thus suffices to verify that for each $s \in S$, a fractional perfect matching $f$ of $\G$ exists for which $f(\{1,s\}) = 1$. In fact, since $A$ is Abelian, the map sending each $a \in A$ to $a^{-1}$ is an automorphism of $\G$, and so it suffices to take just one of $s$ and $s^{-1}$ for each pair of an element $s$ and its inverse from $S$. So pick any $s \in S$. Since $\G$ is not a cycle, there exists an element $s' \in S \setminus \{s, s^{-1}\}$. By Proposition~\ref{pro:3HCext}, there exists a Hamilton cycle $C$ of $\G$ containing the path $P = (s',1,s,ss')$. Deleting the edges of $P$ from $C$ and adding the edge $\{s',ss'\}$ (note that $ss' = s's$ since $A$ is Abelian), we obtain a cycle meeting all the vertices of $\G$ except $1$ and $s$. We can now assign the value of $1/2$ to all the edges of this cycle and $1$ to $\{1,s\}$ to obtain a desired fractional perfect matching of $\G$.

To complete the proof, we only need to show that the graphs from the five families from Proposition~\ref{pro:3HCext} are fractional $1$-extendable. To see that this holds for the examples of even order, we can refer to~\cite{ChaCheYu95}. However, we can also use the result of~\cite{CheQui81} that each edge of a connected Cayley graph of an Abelian group of order at least $3$ is contained in a Hamilton cycle. The graph is thus $1$-extendable (just take every other edge of such a Hamilton cycle). To verify that the graphs $\Circ(2m+1 ; \{\pm 1, \pm 3\})$, $m \geq 3$, are fractional $1$-extendable, let $\G'$ be the subgraph obtained by removing the endvertices of the edge $\{0,1\}$ or $\{0,3\}$. Since the cycle $(2,3,4,\ldots , 2m)$, or $(2,1,4,5,\ldots , 2m)$, respectively, is a spanning subgraph of $\G'$, we can apply Proposition~\ref{pro:kriterij}. Similarly, for $\Circ(m ; \{\pm 1, \pm 2\})$, with $m \geq 7$ odd, removing the endvertices of the edge $\{0,1\}$ or $\{0,2\}$ results in a subgraph, spanned by the $3$-cycle $(4, 5, 6)$ and $(m-5)/2$ independent edges. 
\end{proof}

\subsection{Fractional $2$-extendability - the case of even order}
\label{subsec:even}

The problem of classifying the fractional $2$-extendable Cayley graphs of Abelian groups is considerably more difficult. However, in the case of graphs of even order, the classification of connected $2$-extendable Cayley graphs of Abelian groups of even order from~\cite{ChaCheYu95} is of great help. It states that, besides the obvious examples of cycles of length at least $6$, the only connected Cayley graphs of Abelian groups of even order which are not $2$-extendable, are the examples from items (i), (ii), (iii) and (v) of Proposition~\ref{pro:3HCext}, where $m$ is even in case (v). To classify the fractional $2$-extendable Cayley graphs of Abelian groups of even order we thus only need to see if any of these graphs are fractional $2$-extendable.

\begin{proposition}
\label{pro:2-f-ext-CayAb_even}
Let $\G$ be a connected Cayley graph of an Abelian group of even order at least $6$. Then $\G$ is fractional $2$-extendable if and only if it is not isomorphic to one of the following circulants:
\begin{itemize}
\itemsep = 0pt
\item[{\rm (i)}] $\Circ(2m ; \{\pm 1\})$, $m \geq 3$,
\item[{\rm (ii)}] $\Circ(4m ; \{\pm 1, 2m\})$, $m \geq 2$,
\item[{\rm (iii)}] $\Circ(4m+2 ; \{\pm 2, 2m+1\})$, $m \geq 1$,
\item[{\rm (iv)}] $\Circ(4m+2 ; \{\pm 1, \pm 2m\})$, $m \geq 1$,
\item[{\rm (v)}] $\Circ(2m ; \{\pm 1, \pm 2\})$, $m \geq 3$.
\end{itemize}
\end{proposition}

\begin{proof}
As observed above, we only need to prove that none of the examples from the proposition is fractional $2$-extendable. It is clear that this holds for the cycles $\Circ(2m ; \{\pm 1\})$, $m \geq 3$. Similarly, removing the endvertices of the edges $\{0,1\}$ and $\{3,4\}$ of a graph $\G$ from family (v) results in a graph with an isolated vertex. Proposition~\ref{pro:characterization} now implies that $\G$ is not fractional $2$-extendable.

For each of the graphs $\G$ from families (ii)--(iv) we exhibit a pair of independent edges $e$ and $e'$ of $\G$ such that the subgraph $\G'$, obtained from $\G$ by removing the endvertices of these two edges, has an independent set consisting of more than half of the vertices of $\G'$. Then Corollary~\ref{cor:bip} ensures that $\G$ is not fractional $2$-extendable.

For the examples from family (ii), one can take $e = \{0,1\}$ and $e' = \{2m-1,2m\}$, and observe that the resulting subgraph $\G'$ is bipartite with bipartition parts of different sizes (namely, $2m-1$ and $2m-3$). A similar conclusion can be made by taking $e = \{0,2\}$ and $e' = \{2m-1,2m+1\}$ for the examples from family (iii), and by taking $e = \{0,1\}$ and $e' = \{2m+1,2m+2\}$ for the examples from family (iv).
\end{proof}

\subsection{Fractional $2$-extendability - the case of odd order}\label{ssec: f2e-odd}

We now embark on the project of classifying the fractional $2$-extendable connected Cayley graphs of Abelian groups of odd orders. We omit the trivial case of odd cycles from our considerations. Let $\G = \Cay(A ; S)$ be a connected Cayley graph, where $A$ is an Abelian group of odd order at least $5$ and $|S| > 2$. Suppose that $\G$ is not fractional $2$-extendable. Since $\G$ is vertex-transitive, Proposition~\ref{pro:kriterij} implies existence of disjoint edges $e = \{1,s_1\}$ and $e' = \{a,as_2\}$ for some $a \in A$ and $s_1,s_2 \in S$ such that the subgraph $\G' = \G - \{1,s_1,a,as_2\}$ of $\G$ cannot be covered by a set of pairwise disjoint edges and odd cycles. For further reference we record these assumptions.

\begin{assumption}
\label{ass:CayAb}
Let $A$ be a (multiplicative) Abelian group of odd order at least $5$ and $S \subset A$ be such that $\la S \ra = A$, $S^{-1} = S$, $1 \notin S$ and $|S| \geq 3$. Let $\G = \Cay(A;S)$. We assume that $\G$ is not fractional $2$-extendable and we let $s_1, s_2 \in S$ and $a \in A$ be such that the edges $e = \{1,s_1\}$ and $e' = \{a,as_2\}$ are disjoint and the subgraph $\G' = \G - \{1,s_1,a,as_2\}$ does not have a fractional perfect matching. In other words, it does not have a spanning subgraph consisting of pairwise disjoint edges and odd cycles.
\end{assumption}

In the next series of lemmas, we determine the consequences of Assumption~\ref{ass:CayAb} on the group $A$, the set $S$ and the graph $\G$. As the analysis is rather complex, we divide it into the following three main steps. We first analyze the case that the subgroup $H = \la s_1, s_2\ra$ is a proper subgroup of $A$ (Lemma~\ref{le:CayAb_proper_subgroup}). We then focus on the case that $H = A$, but none of $s_1$ and $s_2$ generates $A$ (Lemma~\ref{le:CayAb_cyclic}). The remaining case when at least one of $s_1$ and $s_2$ generates $A$ (in which case $A$ is cyclic), is the most involved. We therefore divide it into several subcases (Lemma~\ref{le:circ_a=b_1}--Lemma~\ref{le:circ_x_big_a_even_2}). 

Throughout this section, for any $s \in S$ the edges of the form $\{b, bs\}$ will be called {\em $s$-edges}. Moreover, for $s \in S$ and a vertex $b \in A$, the unique cycle through $b$ consisting of $s$-edges will be called the {\em $s$-cycle through} $b$.

\begin{lemma}
\label{le:CayAb_proper_subgroup}
With reference to Assumption~\ref{ass:CayAb}, if $\langle s_1, s_2 \rangle$ is a proper subgroup of the group $A$, then $|A| = 3n$ for some odd integer $n \geq 3$, and one of the following holds:
\begin{itemize}
\itemsep = 0pt
\item[{\rm (i)}] $\G \cong \mathrm{Circ}(3n ; \{\pm 1, \pm 3\})$.
\item[{\rm (ii)}] $\G \cong \mathrm{Cay}(\ZZ_{n} \times \ZZ_3 ; \{\pm (1,0), \pm (1,1)\})$.
\item[{\rm (iii)}] $\G \cong \mathrm{Circ}(3n ; \{\pm 1, \pm (n-1), \pm (n+1)\})$.
\end{itemize}
Moreover, the above graphs are indeed not fractional $2$-extendable.
\end{lemma}

\begin{proof}
Suppose that $H = \langle s_1, s_2\rangle$ is a proper subgroup of $A$ and let $m$ denote its index $[A:H]$ in $A$. The graph $\G$ thus contains a spanning subgraph consisting of $m$ disjoint copies of the connected graph $\G_0 = \Cay(H; H \cap S)$, one for each coset of $H$ in $A$. Note that $\G_0$ is of odd order and (being a connected Cayley graph of an Abelian group) contains a Hamilton cycle. Moreover, the following easy observation will play a crucial role throughout the proof. For any $b \in A$, $s \in S$ and $s' \in S \setminus H$, since $A$ is Abelian, $bss' = bs's$, and so the $s$-edges $\{b,bs\}$ and $\{bs',bs's\}$ are linked by the $s'$-edges $\{b,bs'\}$ and $\{bs,bss'\}$. Observe also that $e$ is an edge of $\G_0$ and that $m \geq 3$ since $m$ is odd. We first prove that $\G_0$ is a cycle and $m = 3$ by a series of three claims.
\medskip

\noindent
{\Claim 1:} $e'$ is not an edge of $\G_0$.\\
Suppose on the contrary that $e'$ is an edge of $\G_0$. We can then find a suitable spanning subgraph of $\G'$ to contradict Assumption~\ref{ass:CayAb} as follows. Since $H$ is a proper subgroup of $A$ and $\la S \ra = A$, there exists an $s \in S \setminus H$. To cover the vertices of $H \cup sH$, we can take the edges $\{s,ss_1\}$, $\{sa,sas_2\}$ and all of the edges of the form $\{h,hs\}$, $h \in H \setminus \{1,s_1,a,as_2\}$ (see the left part of Figure~\ref{fig:merged_example}). For each of the remaining copies of $\G_0$ we can simply take a Hamilton cycle. This contradicts Assumption~\ref{ass:CayAb}, proving our claim.
\begin{figure}[htbp]
\begin{center}
\subfigure{
\begin{tikzpicture}[yscale=0.7, xscale=0.8]
\node[label=270:$H$] at (1.5,0) {};
\node[vtx, fill=black, inner sep = 2pt, label=270:$1$] (L1) at (1, 1) {};
\node[vtx, fill=black, inner sep = 2pt] (L3) at (1, 3) {};
\node[vtx, fill=black, inner sep = 2pt] (L5) at (1, 5) {};
\node[vtx, fill=black, inner sep = 2pt, label=270:$a$] (L7) at (1, 7) {};
\node[vtx, fill=black, inner sep = 2pt] (L9) at (1, 9) {};
\node[vtx, fill=black, inner sep = 2pt, label=270: $s_1$] (R2) at (2, 2) {};
\node[vtx, fill=black, inner sep = 2pt] (R4) at (2, 4) {};
\node[vtx, fill=black, inner sep = 2pt, label=270: $as_2$] (R6) at (2, 6) {};
\node[vtx, fill=black, inner sep = 2pt] (R8) at (2, 8) {};
\draw[opacity = .5, magenta, line width = 3 pt] (L1) -- (R2);
\draw[opacity = .5, magenta, line width = 3 pt] (L7) -- (R6);
\draw[thick] (0,0) -- (3,0) -- (3,10) -- (0,10) -- (0,0);

\node[label=270:$sH$] at (5.5,0) {};
\node[vtx, fill=black, inner sep = 2pt, label=270:$s$] (L1a) at (5, 1) {};
\node[vtx, fill=black, inner sep = 2pt] (L3a) at (5, 3) {};
\node[vtx, fill=black, inner sep = 2pt] (L5a) at (5, 5) {};
\node[vtx, fill=black, inner sep = 2pt, label=270:$sa$] (L7a) at (5, 7) {};
\node[vtx, fill=black, inner sep = 2pt] (L9a) at (5, 9) {};
\node[vtx, fill=black, inner sep = 2pt, label=280: $ss_1$] (R2a) at (6, 2) {};
\node[vtx, fill=black, inner sep = 2pt] (R4a) at (6, 4) {};
\node[vtx, fill=black, inner sep = 2pt, label=270:$sas_2$] (R6a) at (6, 6) {};
\node[vtx, fill=black, inner sep = 2pt] (R8a) at (6, 8) {};

\draw[opacity = .5, blue, line width = 3 pt] (L1a) -- (R2a);
\draw[opacity = .5, blue, line width = 3 pt] (L7a) -- (R6a);
\draw[opacity = .5, blue, line width = 3 pt] (L3) -- (L3a);
\draw[opacity = .5, blue, line width = 3 pt] (L5) -- (L5a);
\draw[opacity = .5, blue, line width = 3 pt] (L9) -- (L9a);
\draw[opacity = .5, blue, line width = 3 pt] (R4) -- (R4a);
\draw[opacity = .5, blue, line width = 3 pt] (R8) -- (R8a);

\draw[thick] (4,0) -- (7,0) -- (7,10) -- (4,10) -- (4,0);
\end{tikzpicture}
}\hspace{3cm}
\subfigure{
\begin{tikzpicture}[yscale=0.8]
\node[label=270:$H$] at (0,0) {};
\node[vtx, fill=blue, inner sep = 3pt, label=180:$s_1^{n-1}$] (08) at (0, 8) {};
\node[vtx, fill=black, inner sep = 2pt, label=180:$s_1^{n-2}$] (07) at (0, 7) {};
\node[vtx, fill=blue, inner sep = 3pt, label=180:$s_1^{n-3}$] (06) at (0, 6) {};
\node[vtx, fill=black, inner sep = 2pt] (05) at (0, 5) {};
\node[vtx, fill=black, inner sep = 2pt] (03) at (0, 3) {};
\node[vtx, fill=blue, inner sep = 3pt, label=180:$s_1^2$] (02) at (0, 2) {};
\node[vtx, fill=black, inner sep = 2pt, label=180:$s_1$] (01) at (0, 1) {};
\node[vtx, fill=black, inner sep = 2pt, label=180:$1$] (00) at (0, 0) {};

\node[label=270:$aH$] at (2,0) {};
\node[vtx, fill=black, inner sep = 2pt, label=45:$ss_1^{n-1}$] (28) at (2, 8) {};
\node[vtx, fill=blue, inner sep = 3pt] (27) at (2, 7) {};
\node[vtx, fill=black, inner sep = 2pt] (26) at (2, 6) {};
\node[vtx, fill=blue, inner sep = 3pt] (25) at (2, 5) {};
\node[vtx, fill=blue, inner sep = 3pt] (23) at (2, 3) {};
\node[vtx, fill=black, inner sep = 2pt, label=45:$ss_1^2$] (22) at (2, 2) {};
\node[vtx, fill=black, inner sep = 2pt, label=45:$ss_1$] (21) at (2, 1) {};
\node[vtx, fill=blue, inner sep = 3pt, label=45:$s$] (20) at (2, 0) {};
\node[label=270:$a^2H$] at (4,0) {};
\node[vtx, fill=blue, inner sep = 3pt, label=0:$s^2s_1^{n-1}$] (48) at (4, 8) {};
\node[vtx, fill=black, inner sep = 2pt] (47) at (4, 7) {};
\node[vtx, fill=blue, inner sep = 3pt] (46) at (4, 6) {};
\node[vtx, fill=black, inner sep = 2pt] (45) at (4, 5) {};
\node[vtx, fill=black, inner sep = 2pt, label=0:$s^2s_1^3$] (43) at (4, 3) {};
\node[vtx, fill=black, inner sep = 2pt] (42) at (4, 2) {};
\node[vtx, fill=blue, inner sep = 3pt, label=0:$s^2s_1$] (41) at (4, 1) {};
\node[vtx, fill=black, inner sep = 2pt, label=0:$s^2$] (40) at (4, 0) {};
\draw[opacity = .5, magenta, line width = 3 pt] (00) -- (01);
\draw[thick] (01) -- (03);
\draw[dashed,thick] (03) -- (05);
\draw[thick] (05) -- (08);

\draw[thick] (20) -- (21);
\draw[opacity = .5, magenta, line width = 3 pt] (21) -- (22);
\draw[thick] (22) -- (23);
\draw[dashed, thick] (23) -- (25);
\draw[thick] (25) -- (28);

\draw[thick] (40) -- (43);
\draw[dashed, thick] (43) -- (45);
\draw[thick] (45) -- (48);

\draw[thick](00) -- (40);
\draw[thick] (01) -- (41);
\draw[thick] (02) -- (42);
\draw[thick] (03) -- (43);

\draw[thick] (05) -- (45);
\draw[thick] (06) -- (46);
\draw[thick] (07) -- (47);
\draw[thick] (08) -- (48);

\draw[thick] (40) to[out=135, in=300] (06);
\draw[thick] (41) to[out=135, in=300] (07);
\draw[thick] (42) to[out=135, in=300] (08);
\draw[thick] (00) -- (43);
\draw[thick] (05) -- (48);

\draw[thick] (00) to[out=100,in=260] (08);
\draw[thick] (20) to[out=100,in=260] (28);
\draw[thick] (40) to[out=100,in=260] (48);
\end{tikzpicture}
}
\caption{Two situations from the proof of Lemma~\ref{le:CayAb_proper_subgroup}.}
\label{fig:merged_example}
\end{center}
\end{figure}
\medskip

\noindent
{\Claim 2:} $aH \neq H$, $S \cap H = \{s_1, s_1^{-1}\}$ and $\G_0$ is a cycle.\\
Since $s_2 \in H$, the edge $e'$ is contained in the copy of $\G_0$ corresponding to the coset $aH$, and so Claim~1 implies that $aH \neq H$. Corollary~\ref{cor:CayAb1fext} ensures that unless $\G_0$ is a cycle of odd length, we can extend each of $e$ and $e'$ to a fractional perfect matching of the corresponding copy of $\G_0$. Since we can take a Hamilton cycle in each of the remaining copies of $\G_0$, this would contradict Assumption~\ref{ass:CayAb}. Therefore, $S \cap H = \{s_1, s_1^{-1}\}$, implying that $\G_0$ is a cycle consisting of $s_1$-edges.
\medskip

\noindent
{\Claim 3:} $m = 3$.\\
Suppose on the contrary that $m > 3$. Note that for each $b \in A$ and each $s \in S \setminus \{s_1, s_1^{-1}\}$, the fact that $s$ is of odd order implies that $s^{-1}bH$, $bH$ and $sbH$ are three different cosets. It is thus easy to see that since $m > 3$, we can find $s,s' \in S \setminus \{s_1, s_1^{-1}\}$ such that $(H \cup sH) \cap (aH \cup s'aH) = \emptyset$. We can then again find a suitable spanning subgraph of $\G'$ to contradict Assumption~\ref{ass:CayAb}. Namely, take the edges $\{s, ss_1\}$ and $\{s'a, s'as_2\}$, cover all of the remaining vertices of $H \cup sH$ (in $\G'$) by $s$-edges, all of the remaining vertices of $aH \cup s'aH$ by $s'$-edges, and take the $s_1$-cycles in each of the remaining $m-4$ copies of $\G_0$.
\medskip

Note that Claim~2 implies that $s_2 \in \{s_1, s_1^{-1}\}$. Since we can replace $a$ by $as_1^{-1}$ if needed, we can thus assume that $s_2 = s_1$. Denote the order of $s_1$ by $n$ and note that Claims~2 and~3 imply that $\G$ is of order $3n$. Moreover, $\G$ consists of the three $n$-cycles corresponding to the $s_1$-edges, which are linked by the $s$-edges corresponding to all $s \in S \setminus \{s_1, s_1^{-1}\}$. We next determine the set $S \setminus \{s_1, s_1^{-1}\}$.
\medskip

\noindent
{\Claim 4:} $S \cap aH \subset \{as_1, as_1^{-1}\}$ and $|S| \in \{4,6\}$.\\
As in the proof of Claim~3 we see that for each $s \in S \setminus \{s_1, s_1^{-1}\}$, one of $s$ and $s^{-1}$ is in $aH$ and the other is in $a^2H$. Take any $s \in S \cap aH$ and let $i \in \{0,1,\ldots , n-1\}$ be such that $a = ss_1^i$. Suppose that $i \notin \{1,n-1\}$. Depending on whether $i$ is odd or even, we can take the edge $\{s_1^{i-1}, ss_1^{i-1}\}$ or $\{s_1^{i+2},ss_1^{i+2}\}$, respectively, cover all of the remaining vertices of $H$ and $aH$ by independent $s_1$-edges, and take the $s_1$-cycle in $a^2H$ to contradict Assumption~\ref{ass:CayAb}. This shows that $S \cap aH \subseteq \{as_1, as_1^{-1}\}$. Since $|S \cap H| = 2$ and $s^{-1} \in S \cap aH$ for each $s \in S \cap a^2H$, this also implies that $|S| \in \{4,6\}$, as claimed. 
\medskip

We now analyze the two possibilities, depending on whether $|S| = 4$ or $|S| = 6$.
\medskip

\noindent
{Case $|S| = 4$:}\\
By Claim~4, $S \cap aH$ consists of one of $as_1$ and $as_1^{-1}$. We only consider the possibility that $s = as_1^{-1} \in S \cap aH$ in all detail (the other one is completely analogous). Let $\ell \in \{0,1,\ldots , n-1\}$ be such that $s^3 = s_1^\ell$. We claim that $\ell \in \{1,n-3,n-1\}$. Suppose on the contrary that $\ell \notin \{1,n-3,n-1\}$ and note that this implies $n > 3$. If $\ell$ is even (in which case $\ell \leq n-5$), consider the odd cycle $C = (s_1^4, ss_1^4, s^2s_1^4, s_1^{\ell + 4}, s_1^{\ell + 3}, \ldots , s_1^5)$. It covers the vertices $ss_1^4$ and $s^2s_1^4$ from $aH$ and $a^2H$, respectively, and all vertices of $H$, except for the endvertices of $e$, the consecutive vertices $s_1^2$ and $s_1^3$, and an even number of consecutive vertices $s_1^{\ell+5}, s_1^{\ell+6}, \ldots, s_1^{-1}$. Since $C$ does not meet $e$ or $e'$, we can thus take the cycle $C$, the edge $\{s^2s_1, s^2s_1^2\}$, and cover all of the remaining vertices of $H$ by independent $s_1$-edges and the remaining vertices of $aH$ and $a^2H$ by $s$-edges to contradict Assumption~\ref{ass:CayAb}. Similarly, if $\ell$ is odd (in which case $\ell \geq 3$), take the cycle $(s_1^{-1},ss_1^{-1},s^2s_1^{-1},s_1^{\ell-1}, s_1^{\ell}, \ldots , s_1^{-2})$, the edge $\{s^2s_1, s^2s_1^2\}$, and cover the remaining vertices by independent $s_1$-edges and $s$-edges, again contradicting Assumption~\ref{ass:CayAb}. Therefore, $\ell \in \{1,n-3,n-1\}$, as claimed.
\smallskip

\noindent
We now show that in this case the graph $\G'$ indeed does not have a fractional perfect matching. To see this take
$$
	U = \{s_1^{2i} \colon 1 \leq i \leq (n-1)/2\} \cup \{ss_1^{2i+1} \colon 1 \leq i \leq (n-3)/2\} \cup \{s^2s_1^{2i} \colon 2 \leq i \leq (n-1)/2\} \cup \{s, s^2s_1^\delta\},
$$
where $\delta = 1$ if $\ell = n-3$ and $\delta = 2$ otherwise (see the right part of Figure~\ref{fig:merged_example}, where the possibility $\ell = n-3$ is depicted). It is not difficult to see that $U$ is an independent set of vertices in $\G'$ containing more than half of its vertices, and so we can apply Corollary~\ref{cor:bip}. To conclude this case we determine the group $A$ and the connection set $S$. If $\ell \in \{1,n-1\}$, then clearly $|s| = 3n$ and $\G \cong \mathrm{Circ}(3n ; \{\pm 1, \pm 3\})$. If however $\ell = n-3$, then $s^3 = s_1^{-3}$ implies that $a^3 = 1$. Since $\langle s_1 \rangle$ and $\langle a \rangle$ clearly intersect trivially, $A \cong \ZZ_n \times \ZZ_3$, where we can let $s_1$ corresponds to $(-1,0)$ and $a$ to $(0,1)$. Since $s = as_1^{-1}$, it is clear that we can assume that $S = \{\pm (1,0), \pm (1,1)\}$. 
\medskip

\noindent
{Case $|S| = 6$:}\\
By Claim~4, $S \cap aH = \{as_1, as_1^{-1}\}$. Denote $s = as_1^{-1}$ and $s' = as_1 = ss_1^{2}$. Again, set $\ell \in \{0,1,\ldots , n-1\}$ be such that $s^3 = s_1^\ell$, and note that the above proof still applies. Thus, $\ell \in \{1,n-3,n-1\}$. 
The claims of the lemma in the case of $n = 3$ can be verified easily, so we assume $n > 3$. We claim that $\ell = n-3$. If $\ell = 1$, then $s^2s_1^{-1}s' = s_1^2$. Consider the odd cycle $C = (s_1^2, s^2s_1^{-1},s^2s_1^{-2}, ss_1^{-2},ss_1^{-1},s_1^{-1},s_1^{-2},\ldots , s_1^3)$. It covers all the vertices of $H$ (except the endvertices of $e$) and the vertices $ss_1^{n-2}, ss_1^{n-1}$ of $aH$ and $s^2s_1^{n-2}, s^2s_1^{n-1}$ of $a^2H$. We can thus take the cycle $C$, the edge $\{s^2s_1, s^2s_1^2\}$, and cover the remaining vertices of $aH$ and $a^2H$ by $s$-edges to contradict Assumption~\ref{ass:CayAb}. Similarly, if $\ell = n-1$, we can take the odd cycle $(s_1^2,ss_1^4,ss_1^5,\ldots, ss_1^{-1},s^2s_1)$, the edges $\{s,s^2s_1^2\}$ and $\{ss_1^3,s^2s_1^5\}$, and cover all of the remaining vertices by independent $s_1$-edges, a contradiction. Therefore, $\ell = n-3$, as claimed. It is not difficult to see that in this case  the above defined set $U$ is still an independent set of vertices of $\G'$. Therefore, Corollary~\ref{cor:bip} implies that $\G'$ has no fractional perfect matching. 
\smallskip

\noindent
Since $\ell = n-3$, the argument from the above case shows that $A \cong \ZZ_n \times \ZZ_3$, where $s_1$ corresponds to $(-1,0)$ and $a$ to $(0,1)$. Note also that $s' = ss_1^2 = as_1$, and so $s'^{-1}$ corresponds to $(1,-1)$. We can thus assume that $S = \{\pm (1,0), \pm (1,1), \pm (1,-1)\}$. That $\mathrm{Cay}(\ZZ_{n} \times \ZZ_3 ; \{\pm (1,0), \pm (1,1), \pm (1,-1)\}) \cong \mathrm{Circ}(3n ; \{\pm 1, \pm (n-1), \pm (n+1)\})$, is easy to verify.
\end{proof}

We remark that the circulant from item (iii) of the above Lemma~\ref{le:CayAb_proper_subgroup} is what is known as the lexicographic product of the $n$-cycle by an edgeless graph of order $3$ (just like the wreath graph from item (iv) of Proposition~\ref{pro:2-f-ext-CayAb_even} is the lexicographic product of the $(2m+1)$-cycle by an edgeless graph of order $2$). Moreover, we mention that if $n$ is coprime to $3$ in item (ii) of Lemma~\ref{le:CayAb_proper_subgroup}, then the graph is a circulant (with connection set being $\{\pm 1, \pm (n-1)\}$ or $\{\pm 1, \pm (n+1)\}$, depending on whether $n \equiv 1 \pmod {3}$ or $n \equiv 2 \pmod{3}$, respectively). 

\begin{lemma}
\label{le:CayAb_cyclic}
With reference to Assumption~\ref{ass:CayAb}, if $A = \langle s_1, s_2\rangle$ but none of $s_1$ and $s_2$ generates $A$, then $A \cong \ZZ_{n} \times \ZZ_3$ for some odd integer $n$ with $3 \mid n$. Moreover, $\G \cong \Cay(\ZZ_{n} \times \ZZ_3 ; \{\pm (1,0), \pm (1,1)\})$ or $\G \cong \Cay(\ZZ_n \times \ZZ_3 ; \{\pm (1,0), \pm (1,1), \pm (1,-1)\}) \cong \mathrm{Circ}(3n ; \{\pm 1, \pm (n-1), \pm (n+1)\})$.
\end{lemma}

\begin{proof}
Suppose that $A = \langle s_1, s_2\rangle$ but none of $s_1$ and $s_2$ generates $A$. Denote the order of $s_1$ by $n_1$ and the index $[A : \langle s_1 \rangle]$ by $m_1$, and similarly $n_2,m_2$ for $s_2$. Then $A$ (and thus $\G$) is of order $m_1 n_1 = m_2 n_2$. Since none of $s_1, s_2$ generates $A$ and  $m_1,m_2$ are odd, we also have $m_1,m_2\geq 3$.  Note that the edge $e'=\{a,as_2\}$ connects vertices from two different cosets  $a\langle s_1 \rangle, as_2\langle s_1\rangle$ in $A$, as equality would imply $s_2\in\langle s_1\rangle= \langle s_1,s_2\rangle$, a contradiction. 
Since $\langle s_1 \rangle$ and $\langle s_1, s_2\rangle$ are of orders $n_1$ and $m_1 n_1$, respectively, there exists a unique $\ell \in \{0,1,\ldots , n_1 - 1\}$ such that $s_2^{m_1} = s_1^\ell$. The order $n_2$ of $s_2$ thus equals $m_1 n_1/\gcd(\ell, n_1)$. Since $m_1n_1 = m_2n_2$, it thus follows that $\gcd(\ell, n_1) = m_2$ (implying in particular that $m_2$ divides $n_1$). 
We proceed to show that $m_1 = m_2 = 3$, that $3 \mid n_1$, and that $a \in (s_2\langle s_1 \rangle) \cap (s_1^2\langle s_2 \rangle)$ by the following three claims.
\medskip

\noindent
{\Claim 1:} $a, as_2 \notin \langle s_1 \rangle$.\\
Suppose first that $a \in \langle s_1 \rangle$. Than we can simply take the edge $\{s_2, s_2s_1\}$ (which is of course disjoint from $e'$ since $e$ is), all of the edges $\{s_1^i, s_1^is_2\}$, where $s_1^i \notin \{1,s_1,a\}$, together with the $s_1$-cycles (of odd length $n_1$) on all of the remaining $m_1 - 2$ cosets of $\langle s_1 \rangle$ to contradict Assumption~\ref{ass:CayAb}. The case when $as_2 \in \langle s_1 \rangle$ is similar.
\medskip

\noindent
{\Claim 2:} $m_1=m_2=3$ and $3 \mid n_1$.\\
First, consider the $s_2$-cycle $C$ through $s_1^2$. Since $\gcd(\ell,n_1) = m_2$, the set of vertices of some coset $s_2^i\langle s_1 \rangle$, where $0 \leq i < m_1$, that are on $C$, is the set $\{s_2^i s_1^{2+jm_2} \colon 0 \leq j < n_1/m_2\}$. Note that $m_2 \geq 3$ implies that the cycle $C$ is disjoint from $e$. Moreover, the set of elements of a coset of $\la s_1 \ra$, not covered by $C$, consists of $n_1/m_2$ sets of $m_2-1$ consecutive vertices. Therefore, if $C$ does not contain $e'$, we can take the cycle $C$ (which is of odd length $n_2$) and all of the remaining $s_2$-edges between the vertices of $a\langle s_1 \rangle$ and $as_2\langle s_1 \rangle$. As for the cosets of $\langle s_1 \rangle$ in $A$ except $a\langle s_1 \rangle$ and $as_2\langle s_1 \rangle$, we can cover their remaining vertices (those not already on $C$) with disjoint $s_1$-edges (since $m_2$ is odd). This contradicts Assumption~\ref{ass:CayAb}, showing that $C$ does contain $e'$. In a s similar way we see that the $s_2$-cycle $C'$ containing $s_1^{-1}$ also contains $e'$, and thus $C' = C$. This implies that $m_2 = 3$. Exchanging the roles of $s_1$ and $s_2$, we have $m_1 = m_2 = 3$. Since $m_2$ divides $n_1$, the claim follows.
\medskip

\noindent
{\Claim 3:} $a \in (s_2\langle s_1 \rangle) \cap (s_1^2\langle s_2 \rangle)$.\\
Since $e'$ is contained in $C$, it clearly follows that $a \in s_1^2\langle s_2 \rangle$. Moreover, since $m_2 = 3$ and $e'$ is disjoint from $\la s_1 \ra$, we also get $a \in s_2\langle s_1 \rangle$, which proves our claim.
\smallskip

\noindent
That the conclusion of the lemma is correct when $n_1 = 3$ (in which case $\ell = 0$ and $A \cong \ZZ_3 \times \ZZ_3$) can easily be verified. Let us thus assume that $n_1 > 3$ and note that in this case $n_1 \geq 9$ (since $m_2 \mid n_1$). Recall that $\gcd(\ell, n_1) = 3$. It can be verified that renaming the vertices using multiplication by $s_1^{-1}$ and then exchanging the roles of $s_1$ and $s_1^{-1}$, the parity of the corresponding $\ell$ changes. This shows that we can assume that $\ell$ is odd, which implies that $3 \leq \ell \leq n_1 - 6$ (recall that $n_1$ is odd). We now prove that $a = s_1^{-1}s_2$, that $\ell = 3$ and that $A \cong \ZZ_{n_1} \times \ZZ_3$ by the following two claims.
\medskip

\noindent
{\Claim 4:} $a = s_1^{-1}s_2$.\\
Consider the cycle $C = (s_1^{-1}, s_1^{-1}s_2, s_1^{-1}s_2^2, s_1^{\ell-1}, s_1^{\ell}, \ldots , s_1^{-2})$. It is of odd length $n-\ell+3$ and covers all of the vertices of $\la s_1 \ra$, except for the endvertices of $e$ and the $\ell - 3$ consecutive vertices $s_1^2, s_1^3, \ldots , s_1^{\ell-2}$. Therefore, if $C$ does not contain $e'$, then we can take this cycle, cover the remaining vertices of $\langle s_1 \rangle$ by independent $s_1$-edges and cover the remaining vertices of $s_2\langle s_1 \rangle$ and $s_2^2\langle s_1 \rangle$ by $s_2$-edges to contradict Assumption~\ref{ass:CayAb}. It thus follows that $C$ contains $e'$, forcing $a = s_1^{-1}s_2$, as claimed.
\medskip

\noindent
{\Claim 5:} $\ell = 3$ and $A \cong \ZZ_{n_1} \times \ZZ_3$.\\
Suppose on the contrary that $\ell \neq 3$, consider the cycle $(s_1^{-3}, s_1^{-3}s_2, s_1^{-3}s_2^2, s_1^{\ell-3}, s_1^{\ell-2}, \ldots , s_1^{-4})$, and note that it does not meet $e$ or $e'$. Similarly as above, we see that we can cover the remaining vertices of $\G'$ by independent $s_1$-edges and $s_2$-edges to contradict Assumption~\ref{ass:CayAb}. Therefore, $\ell = 3$, and so $s_1^{-1}s_2$ is of order $3$. Since $\langle s_1 \rangle$ and $\langle s_1^{-1}s_2 \rangle$ intersect trivially, $A \cong \ZZ_{n_1} \times \ZZ_3$, where $s_1$ and $s_2$ correspond to $(1,0)$ and $(1,1)$, respectively.
\medskip

If $|S| = 4$, the proof is complete. Suppose then that $|S| > 4$. If $S$ contains some $s_1^i$, $1 < i < n_1-1$, then the subgraph of $\G$ induced on $\langle s_1 \rangle$ is fractional $1$-extendable by Corollary~\ref{cor:CayAb1fext}. However, this contradicts Assumption~\ref{ass:CayAb}, since we can cover the vertices of $s_2\langle s_1 \rangle$ and $s_2^2\langle s_1 \rangle$ in $\G'$ by $s_2$-edges. We thus have an element $s = s_1^i s_2 \in S$ with $1 \leq i < n_1$. We claim that $i = n_1 - 2$. Suppose on the contrary that this is not the case. If $i = n_1-1$, we can take the $7$-cycle $(s_1^{-1}, s_1^{-2}s_2, s_1^{-2}s_2^2, s_1^{-3}s_2^2, s_1^{-3}s_2, s_1^{-3}, s_1^{-2})$ and cover the remaining vertices of $\G'$ by independent $s_1$-edges, contradicting Assumption~\ref{ass:CayAb}. If $i$ is even with $i < n_1-1$, we can take the cycle $(s_1^{-3}, s_1^{-3}s_2, s_1^{-3}s_2^2, s_1^i, s_1^{i+1}, \ldots , s_1^{-4})$, cover the remaining vertices of $\langle s_1 \rangle$ by independent $s_1$-edges (we can do this since $i$ is even), and cover the remaining vertices of $s_2\langle s_1 \rangle$ and $s_2^2\langle s_1 \rangle$ by $s_2$-edges, a contradiction. Finally, if $i$ is odd with $i < n_1 - 2$, we can take the cycle $(s_1^2, s_1^{i+2}s_2, s_1^{i+1}s_2, \ldots , s_1^2 s_2)$ and cover all of the remaining vertices by independent $s_1$-edges, a contradiction. This proves our claim that $i = n_1 - 2$. Note that this actually forces $S \cap s_2\la s_1\ra = \{s_1^{-2}s_2\}$, and so $S = \{s_1, s_1^{-1}, s_2, s_2^{-1}, s_1^{-2}s_2, s_1^{-1}s_2^2\}$. Of course, the element $s_1^{-1}s_2^2$ corresponds to $(1,-1)$ in $\ZZ_{n_1} \times \ZZ_3$.
\end{proof}

To conclude the analysis of fractional $2$-extendability of Cayley graphs of Abelian groups, we finally consider the examples in which at least one of $s_1$ and $s_2$ from Assumption~\ref{ass:CayAb} generates $A$. In particular, $\G$ is a circulant in this case. We can thus assume that $A = \ZZ_n$ for some odd integer $n \geq 5$ (with additive operation), and that one of $s_1$ and $s_2$, say $s_2$, is $1$. It is easy to see that none of the graphs $\Circ(n ; \{\pm 1, \pm 2\})$ and $\Circ(n ; \{\pm 1, \pm 3\})$ is fractional $2$-extendable. For the first one we can repeat the argument from the proof of Proposition~\ref{pro:2-f-ext-CayAb_even}. For the second one, removing the endvertices of the edges $\{0,1\}$ and $\{2,3\}$ results in a bipartite graph of odd order, and we can apply Corollary~\ref{cor:bip}. We shall thus only consider the examples $\Circ(n ; S)$, where $1 \in S$ but either $S \cap \{2,3\} = \emptyset$ or $|S| > 4$.

It turns out that the analysis of the situation is somewhat different if $s_1 \in \{\pm 1\}$ or $s_1 \notin \{\pm 1\}$. We first consider the former case. Observe that, replacing $e' = \{a, a+1\}$ by $\{a+1, a+2\}$ if necessary, we can assume that $s_1 = 1$. Moreover, exchanging the roles of $e$ and $e'$ if necessary, we can assume that $a \leq (n-1)/2$. In the following two lemmas we will thus be working with the following assumption.

\begin{assumption}
\label{ass:circ1}
Let $n \geq 5$ be an odd integer, let $S \subset \ZZ_n \setminus \{0\}$ be such that $-S = S$, $1 \in S$ and that $|S| \geq 4$, and that either $S \cap \{2,3\} = \emptyset$ or $|S| \geq 6$. Let $\G = \Circ(n;S)$. We assume that $a \in \ZZ_n$, where $2 \leq a \leq (n-1)/2$, is such that the subgraph $\G' = \G - \{0,1,a,a+1\}$ does not have a fractional perfect matching. We also denote $e = \{0,1\}$ and $e' = \{a,a+1\}$. 
\end{assumption}

\begin{lemma}
\label{le:circ_a=b_1}
With reference to Assumption~\ref{ass:circ1}, $a$ is odd and there exists no odd $s \in S$ with $1 < s \leq (n-1)/2$.
\end{lemma}

\begin{proof}
We first show that there is no odd $s \in S$ with $1 < s \leq (n-1)/2$. Suppose on the contrary that such an $s$ exists. If $a$ is odd, we let $j \in \{2,3,\ldots , a-1\}$ be the smallest even integer such that $j+s \geq a+2$. It is easy to see that in this case $j + s \leq n+2-s$. Consider the cycle $C = (2,3,\ldots , j, j+s,j+s+1,\ldots , n+2-s)$ of $\G'$ and note that it is of odd length $n-2s+2$ (see the left part of Figure~\ref{fig:lemma_3_8}, where the case of $n = 21$, $a = 9$ and $s = 7$ is presented). The vertices of $\G'$, which are not covered by $C$, are the $a-j-1$ consecutive vertices $j+1, j+2, \ldots, a-1$, the $j+s-a-2$ consecutive vertices $a+2, a+3, \ldots, j+s-1$, and the $s-3$ consecutive vertices $n+3-s, n+4-s, \ldots, n-1$. Since all of these three numbers are even, we can cover these vertices by independent $1$-edges, contradicting Assumption~\ref{ass:circ1} by Proposition~\ref{pro:kriterij}. 
\begin{figure}[htbp]
\begin{center}
\subfigure{
\begin{tikzpicture}[scale = 0.6]
\foreach \j in {0,1,2,3,4,5,6,7,8,9,10,11,12,13,14,15,16,17,18,19,20}{
\node[vtx, fill=black, inner sep = 2pt] (A\j) at (360*\j/21 + 270: 4) {};
}
\node[label=270:{$0$}] (L0) at (270: 4.1) {};
\node[label=280:{$1$}] (L1) at (270 + 360/21: 4.1) {};
\node[label=280:{$e$}] (Le) at (274: 2.9) {};
\node[label=280:{$2$}] (L2) at (270 + 360*2/21: 4.1) {};
\node[label=0:{$j$}] (L3) at (270 + 360*4/21: 4.1) {};
\node[label=10:{$a$}] (L4) at (270 + 360*9/21: 4) {};
\node[label=80:{$a+1$}] (L5) at (270 + 360*10/21: 4) {};
\node[label=20:{$e'$}] (Lee) at (80: 2.7) {};
\node[label=92:{$j+s$}] (L6) at (270 + 360*11/21: 3.9) {};
\node[label=180:{$n+2-s$}] (L7) at (270 + 360*16/21: 4.1) {};
\node[label=270:{$n-1$}] (L8) at (270 + 360*20/21: 4) {};
\node[label=0:{$n = 21$}] (n) at (180: 1.1) {};
\draw[magenta, line width = 3 pt] (A0) -- (A1);
\draw[magenta, line width = 3 pt] (A9) -- (A10);
\draw[blue, line width = 2 pt] (A2) -- (A3) -- (A4) -- (A11) -- (A12) -- (A13) -- (A14) -- (A15) -- (A16) -- (A2);
\draw[blue, line width = 2 pt] (A5) -- (A6);
\draw[blue, line width = 2 pt] (A7) -- (A8);
\draw[blue, line width = 2 pt] (A17) -- (A18);
\draw[blue, line width = 2 pt] (A19) -- (A20);
\end{tikzpicture}
}
\hspace{8mm}
\subfigure{
\begin{tikzpicture}[scale = 0.6]
\foreach \j in {0,1,2,3,4,5,6,7,8,9,10,11,12,13,14,15,16,17,18}{
\node[vtx, fill=black, inner sep = 2pt] (A\j) at (360*\j/19 + 270: 4) {};
}
\node[label=270:{$0$}] (L0) at (270: 4.1) {};
\node[label=280:{$1$}] (L1) at (270 + 360/19: 4.1) {};
\node[label=280:{$e$}] (Le) at (274: 2.9) {};
\node[label=280:{$2$}] (L2) at (270 + 360*2/19: 4.1) {};
\node[label=10:{$a$}] (L4) at (270 + 360*8/19: 4) {};
\node[label=80:{$a+1$}] (L5) at (270 + 360*9/19: 4) {};
\node[label=20:{$e'$}] (Lee) at (80: 2.7) {};
\node[label=92:{$a+2$}] (L6) at (270 + 360*10/19: 4) {};
\node[label=180:{$a+s+2$}] (L7) at (270 + 360*14/19: 4.1) {};
\node[label=270:{$n-1$}] (L8) at (270 + 360*18/19: 4) {};
\node[label=0:{$n = 19$}] (n) at (180: 1.1) {};
\draw[magenta, line width = 3 pt] (A0) -- (A1);
\draw[magenta, line width = 3 pt] (A8) -- (A9);
\draw[blue, line width = 2 pt] (A10) -- (A11) -- (A12) -- (A13) -- (A14) -- (A10);
\draw[blue, line width = 2 pt] (A2) -- (A3);
\draw[blue, line width = 2 pt] (A4) -- (A5);
\draw[blue, line width = 2 pt] (A6) -- (A7);
\draw[blue, line width = 2 pt] (A15) -- (A16);
\draw[blue, line width = 2 pt] (A17) -- (A18);
\end{tikzpicture}
}
\caption{Two situations from the proof of Lemma~\ref{le:circ_a=b_1}.}
\label{fig:lemma_3_8}
\end{center}
\end{figure}

Similarly, if $a$ is even and $a \neq 2$, then there exists an odd $j \in \{3,4,\ldots , a-1\}$ such that $a+2 \leq j+s \leq n+2-s$. We can then again take the cycle $(2,3,\ldots , j, j+s,j+s+1,\ldots , n+2-s)$ and cover the remaining vertices of $\G'$ by independent $1$-edges. For the case when $a = 2$, we can do the following. We take an even $s \in S$ with $s \leq n-5$, which exists since otherwise $S = \{\pm 1, \pm 3\}$, contradicting Assumption~\ref{ass:circ1}. We can then take the cycle $(4,5,\ldots , 4+s)$ and cover the remaining vertices of $\G'$ by independent $1$-edges. This finally shows that there is no odd $s \in S$ with $1 < s \leq (n-1)/2$, as claimed.

To complete the proof, we need to show that $a$ is odd. Suppose on the contrary that $a$ is even. Then $\G'$ has an even number of consecutive vertices between $2$ and $a-1$. Since $s \geq (n+1)/2$ for each odd $s \in S$ with $s \neq 1$, it follows that $s \leq (n-1)/2$ for each even $s \in S$ with $s \neq n-1$. If there exists an even $s \in S$ with $s \leq n-a-3$, then we can take the odd cycle $(a+2,a+3,\ldots , a+s+2)$ in $\G'$, and cover the remaining vertices of $\G'$ by independent $1$-edges, contradicting Assumption~\ref{ass:circ1} (see the right part of Figure~\ref{fig:lemma_3_8}, where the case of $n = 19$, $a = 8$ and $s = 4$ is presented). Thus $s > n-a-3$ for each even $s \in S$. Since $n-a-3$ is even, this in fact implies that $s \geq n-a-1 \geq (n-1)/2$ (recall that $a \leq (n-1)/2$) for each even $s \in S$. But then $S = \{\pm 1, \pm (n-1)/2\}$ and $a = (n-1)/2$, and so $(n-1, a-1, a-2,\ldots , 2, a+2, a+3, \ldots, n-2)$ is a Hamilton cycle of $\G'$, a contradiction.
\end{proof}

\begin{lemma}
\label{le:circ_a=b_2}
With reference to Assumption~\ref{ass:circ1}, $\G$ is isomorphic to $\Circ(n ; \{\pm 1, \pm 3\})$, or $n = 3m$ for some odd $m \geq 3$, and $\G$ is isomorphic to one of the following graphs:
\begin{itemize}
\itemsep = 0pt
\item[{\rm (i)}] $\Circ(n ; \{\pm 1, \pm (m-1)\})$, 
\item[{\rm (ii)}] $\Circ(n ; \{\pm 1, \pm (m+1)\})$,
\item[{\rm (iii)}] $\Circ(n ; \{\pm 1, \pm (m-1), \pm (m+1)\})$.
\end{itemize}
Moreover, none of these graphs is fractional $2$-extendable.
\end{lemma}

\begin{proof}
Lemma~\ref{le:circ_a=b_1} implies that $a$ is odd (and thus $a \geq 3$) and that there exists no odd $s \in S$ with $1 < s \leq (n-1)/2$. Pick any $s \in S \setminus \{1\}$ with $s \leq (n-1)/2$, which thus must be even. We proceed by proving two claims.
\medskip

\noindent
{\Claim 1:} $s\in\{a-1,a+1\}$.\\
Suppose on the contrary that this does not hold. If $s < a - 1$ then $s \leq a - 3$, and so $(2,3,\ldots , s+2)$ is a an odd cycle in $\G'$. Since we can clearly cover the remaining vertices of $\G'$ by independent $1$-edges, this contradicts Assumption~\ref{ass:circ1}. Thus $s \geq a - 1$. Furthermore, if $s > a + 1$, then $s \geq a + 3$, and so the fact that $s \leq (n-1)/2$ implies that 
$$
	a + 2 \leq s - 1 < s + 2 < n - s + 2 < n + a - s + 2 \leq n - 1.
$$
We can thus take the cycle $(2,s+2,s+3, \ldots , n-s+2)$ of odd length $n-2s+2$, every other edge of the cycle $(n-1,s-1,s-2,\ldots, a+2,n+a-s+2,n+a-s+3,\ldots, n-2)$ of even length $2(s-a-2)$, and cover all of the remaining vertices of $\G'$ by independent $1$-edges, contradicting Assumption~\ref{ass:circ1}. Thus $s \in \{a-1,a+1\}$, as claimed. 
\medskip

\noindent
{Claim 2:} $2s \in \{n - a - 2, n - a, n - a + 2\}$ or $s = 2$.\\
If $2s \geq n-a+4$, then $a+s-1 > n+2-s$. The fact that $s+a < n-1$ (and thus $n+2-s > a+3$) then implies that $(2,3,\ldots , a-1,a+s-1,a+s-2,\ldots , n+2-s)$ is a cycle of odd length $2s+2a-n-4 \geq 3$ in $\G'$. Since this leaves us with $n-s-a$ (an even number) consecutive vertices $a+2, a+3, \ldots, n+1-s$ and the same number of consecutive vertices $a+s, a+s+1, \ldots, n-1$, this contradicts Assumption~\ref{ass:circ1}. Thus $2s \leq n-a+2$. Next, if $2s \leq n-a-4$ and $s \neq 2$, then 
$$
	a+2 < a+s-1 < a+s+2 < n-s-1 < n-s+2 < n-1,
$$
and so $(2,3,\ldots , a-1, a+s-1,a+s-2,\ldots, a+2,a+s+2,a+s+3,\ldots , n-s-1, n-1,n-2,\ldots , n-s+2)$ is a cycle of $\G'$, containing all the vertices of $\G'$ except for the endvertices of the two edges $\{a+s, a+s+1\}$ and $\{n-s,n-s+1\}$. Since this again contradicts Assumption~\ref{ass:circ1}, our claim is proven. 
\medskip

We are now ready to determine the set $S$. Suppose first that $2 \in S$. Then Claim~1 implies that $a = 3$. Since Assumption~\ref{ass:circ1} guarantees that $|S| \geq 6$, we have an even $s \in S$ with $2 < s \leq (n-1)/2$. But then Claim~1 implies that $s = 4$, while Claim~2 implies that $8 \in \{n-5,n-3,n-1\}$, forcing $n \in \{9,11,13\}$. In the case that $n \in \{11,13\}$, we have a Hamilton cycle $(2,6,5,7,8,\ldots , n-3,n-1,n-2)$ of $\G'$, a contradiction. Thus, $n = 9$ and $S = \{\pm 1, \pm 2, \pm 4\}$. By Lemma~\ref{le:CayAb_proper_subgroup}, the graph $\Circ(9 ; \{\pm 1, \pm 2, \pm 4\}$ is indeed not fractional $2$-extendable.

We are left with the possibility that $2 \notin S$. Claim~2 implies that $|S| \leq 6$. Suppose first that $|S| = 6$. Then $S = \{\pm 1, \pm (a-1), \pm (a+1)\}$, and so Claim~2 implies that $2(a-1), 2(a+1) \in  \{n-a-2,n-a,n-a+2\}$. Therefore, $2(a-1) = n-a-2$, and thus $n = 3a$. Consequently, $\G \cong \Circ(3a ; \{\pm 1, \pm (a-1), \pm (a+1)\})$, which by Lemma~\ref{le:CayAb_proper_subgroup} is indeed not fractional $2$-extendable. 

Suppose finally, that $|S| = 4$, and let $s \in S$ be the unique element with $1 < s \leq (n-1)/2$. Then Claims~1 and~2 imply that either $s = a - 1$ and $n \in \{3a - 4, 3a - 2, 3a\}$, or $s = a+1$ and $n \in \{3a, 3a + 2, 3a + 4\}$. If $n = 3a$, then $S$ is one of $\{\pm 1, \pm (a-1)\}$ and $\{\pm 1, \pm (a+1)\}$. In any of these two cases, consider the set
$$
	U = \{2,4,\ldots, a-1\} \cup \{a+2,a+4,\ldots , 2a-1\} \cup \{2a+2, 2a + 4, \ldots , 3a-1\}.
$$ 
One can verify that $U$ is an independent set of vertices in $\G'$ and that it contains more than half of the vertices of $\G'$. Corollary~\ref{cor:bip} thus implies that $\G$ is indeed not fractional $2$-extendable.
If however, $n$ is coprime to $3$, then multiplication by $3$ is an isomorphism from $\G$ to $\Circ(n ; 3S)$, where $3S = \{3s \colon s \in S\}$. It is easy to see that in each of the four cases $3S = \{\pm 1, \pm 3\}$. Therefore, $\G \cong \Circ(n ; \{\pm 1, \pm 3\})$, which we already know is not fractional $2$-extendable. 
\end{proof}

We finally analyze the situation in which $s_1 \notin \{\pm 1\}$. Observe that the permutation, mapping each $i \in \ZZ_n$ to $-i$, is an automorphism of a circulant $\Circ(n;S)$. Changing $a$ if necessary, we can thus assume that $s_1 \leq (n-1)/2$. Our assumption throughout the rest of this section will thus be as follows.

\begin{assumption}
\label{ass:circ2}
Let $n \geq 5$ be an odd integer, let $s_1 \in \ZZ_n$ be such that $1 < s_1 \leq (n-1)/2$, and let $S \subset \ZZ_n \setminus \{0\}$ be such that $-S = S$, that $1, s_1 \in S$, and that either $S \cap \{2,3\} = \emptyset$ or $|S| \geq 6$. Let $\G = \Circ(n ; S)$. We assume that $a \in \ZZ_n \setminus \{n-1,0,s_1-1,s_1\}$ is such that the subgraph $\G' = \G - \{0,s_1,a,a+1\}$ does not have a fractional perfect matching. We also denote $e = \{0,s_1\}$ and $e' = \{a, a+1\}$. 
\end{assumption}

We point out that, unlike in Assumption~\ref{ass:circ1}, we now make no assumption on whether $a \leq (n-1)/2$ or not. Just like in all of the previous cases, we now determine the consequences of Assumption~\ref{ass:circ2}, in this case on the parameters $n$, $s_1$ and $a$. In the next two lemmas we first prove that $a > s_1$ and that $s_1$ is even. We then complete the analysis in Lemmas~\ref{le:circ_x_big_a_even} and~\ref{le:circ_x_big_a_even_2}.

\begin{lemma}
\label{le:circ_ell_small}
With reference to Assumption~\ref{ass:circ2}, $a > s_1$.
\end{lemma}

\begin{proof}
Suppose on the contrary that $a < s_1$. Let $V_1 = \{1,2,\ldots , a-1\}$, $V_2 = \{a+2, a+3, \ldots , s_1-1\}$ and $V_3 = \{s_1+1, s_1+2, \ldots , n-1\}$. Observe that renaming the vertices via the map $i \mapsto s_1 - i$ exchanges the names of the endvertices of $e$, as well as the sets $V_1$ and $V_2$. We will refer to this renaming as ``exchanging the roles of $0$ and $s_1$''. We consider two cases depending on the parity of $s_1$.
\medskip

\noindent
{Case} $s_1$ is even:\\
In this case $s_1 \geq 4$ and $|V_3|$ is even, while precisely one of $|V_1|$ and $|V_2|$ is odd. Exchanging the roles of $0$ and $s_1$ if necessary, we can assume that $|V_2|$ is odd. Note that this implies that $a$ is odd and that $V_2$ contains at least the vertex $s_1-1$. Since $2s_1 \leq n-1$, we can take the odd cycle $(n-1, s_1-1,2s_1-1,2s_1, \ldots , n-2)$, which covers the ``last'' vertex of $V_2$ and all of the vertices of $V_3$, except the ``first'' $s_1 - 2$ of them. Since this is an even number, we can clearly cover all of the remaining vertices of $\G'$ by independent $1$-edges, contradicting Assumption~\ref{ass:circ2}.
\medskip

\noindent
{Case} $s_1$ is odd:\\
If $s_1 = 3$ (in which case $a = 1$ is forced), then Assumption~\ref{ass:circ2} guarantees that $S$ contains some even $s$ with $s < n - 4$. We can then take the cycle $(4,5,\ldots , s+4)$ of odd length $s+1$, and cover the remaining vertices of $\G'$ by independent $1$-edges, a contradiction. Therefore, $s_1 \geq 5$. Exchanging the roles of $0$ and $s_1$ if necessary, we can assume that $V_2 \neq \emptyset$. Since this time $|V_3|$ is odd, $|V_1|$ and $|V_2|$ are of the same parity. If $|V_2|$ is even, we can take the cycle $(n-1,s_1-1,s_1-2,2s_1-2,2s_1-1,\ldots , n-2)$, which covers the ``last'' two vertices of $V_2$ and all but the ``first'' $s_1-3$ of $V_3$. Since this is an even number, we can cover the remaining vertices of $\G'$ by independent $1$-edges, a contradiction. If however $|V_2|$ is odd, we take the odd cycle $(n-1,s_1-1,2s_1-1,2s_1,\ldots, n-2)$ and the edge $\{1,s_1+1\}$, which together cover the ``first'' and the ``last'' vertex of $V_1$ and $V_2$, respectively, and all but the $s_1-3$ consecutive vertices $s_1+2, s_1+3, \ldots, 2s_1-2$ of $V_3$. We can thus cover the remaining vertices of $\G'$ by independent $1$-edges, a contradiction.
\end{proof}

\begin{lemma}
\label{le:circ_x_big_a_odd}
With reference to Assumption~\ref{ass:circ2}, $s_1$ is even.
\end{lemma}

\begin{proof}
By way of contradiction, suppose that $s_1$ is odd. By Lemma~\ref{le:circ_ell_small}, $a > s_1$. Let $V_1 = \{1,2,\ldots , s_1-1\}$, $V_2 = \{s_1+1, s_1+2, \ldots , a-1\}$ and $V_3 = \{a+2, a+3, \ldots , n-1\}$.  Then $|V_1|$ is even and precisely one of $|V_2|$ and $|V_3|$ is odd. Exchanging the roles of the vertices $0$ and $s_1$ (which this time exchanges the roles of the sets $V_2$ and $V_3$) if necessary, we can assume that $|V_2|$ is odd. Then $a$ is also odd. We first prove the following claim.
\medskip

\noindent
{Claim:} At least one of $n = 2s_1+1$ and $V_3 = \emptyset$ holds.\\
Suppose that there exists an even integer $j$ with $s_1 < j < a$ such that $a+2 \leq s_1+j < n$. We can then take the cycle $C = (1, s_1+1, s_1+2, \ldots , j, s_1+j, s_1+j+1, \ldots , n-1, s_1-1, s_1-2, \ldots , 2)$. The only vertices of $\G'$ which are not covered by $C$, are the ``last'' $a-j-1$ (which is an even number) vertices of $V_2$ and the ``first'' $s_1+j-a-2$ (also an even number) vertices of $V_3$. These can be covered by independent $1$-edges, contradicting Assumption~\ref{ass:circ2}. Therefore, for each even $j$ with $s_1 < j < a$, at least one of $s_1+j < a+2$ and $s_1+j \geq n$ holds. Since $s_1+1$ is even and $2s_1+1 \leq n$, this proves our claim.
\medskip

To finish the proof, we first consider the possibility that $s_1 > 3$. If $V_3 = \emptyset$ (which occurs if and only if $a = n-2$), we can take the cycle $(1,s_1+1,s_1+2,\ldots , n-3, s_1-3, s_1-4,\ldots, 2)$ and the edge $\{s_1-2,s_1-1\}$ to cover all the vertices of $\G'$, a contradiction. Hence $V_3 \neq \emptyset$, and so $n = 2s_1+1$. If $a \neq s_1+2$, we can take the $5$-cycle $(1,s_1+1,s_1+2,s_1+3,2)$ and cover all of the remaining vertices of $\G'$ by independent $1$-edges, a contradiction. Therefore, $a = s_1+2$. But in this case we can take the $3$-cycle $(n-1,s_1-1,n-2)$ and the edge $\{1,s_1+1\}$, which leaves us with $s_1-3$ consecutive vertices of $V_1$ and $s_1-5$ consecutive vertices of $V_3$. Since these are even numbers, we can cover these vertices by independent $1$-edges, a contradiction.

This leaves us with the possibility that $s_1 = 3$. Note that $V_3 = \emptyset$ if and only if $a = n-2$. Moreover, if $n = 2s_1+1 = 7$, we must have that $a = 5$, again implying that $a = n-2$. By Assumption~\ref{ass:circ2} there exists some even $s \in S$ with $2 \leq s \leq n-5$. We can thus take the cycle $(1,4,5,\ldots , s+2,2)$ and cover the remaining consecutive $n-s-5$ (which is an even number) vertices of $V_2$ by independent $1$-edges. This contradicts Assumption~\ref{ass:circ2}, thus proving that $s_1$ is even, as claimed.
\end{proof}

The above two lemmas thus imply that $s_1$ is even and $a > s_1$. We now show that $s_1$ is roughly $n/3$ and $a$ is roughly $2n/3$.

\begin{lemma}
\label{le:circ_x_big_a_even}
With reference to Assumption~\ref{ass:circ2}, one of the following holds:
\begin{itemize}
\itemsep = 0pt
	\item $n = 3s_1-3$ and $a = 2s_1-2$,
	\item $n = 3s_1-1$ and $a \in \{2s_1-2, 2s_1-1, 2s_1\}$,
	\item $n = 3s_1+1$ and $a \in \{2s_1-1, 2s_1, 2s_1+1\}$,
	\item $n = 3s_1+3$ and $a = 2s_1+1$.
\end{itemize}	
\end{lemma}

\begin{proof}
Recall that $s_1$ is even and $a > s_1$. We let $V_1$, $V_2$ and $V_3$ be as in the proof of Lemma~\ref{le:circ_x_big_a_odd}. Note that this time $|V_1|$ is odd and $|V_2|$ and $|V_3|$ are of the same parity. Exchanging the roles of the vertices $0$ and $s_1$ if necessary, we can assume that $|V_2| \geq |V_3|$ (but see the last paragraph of the proof). Therefore, $a-s_1 \geq n-a-1$, which is equivalent to $2a - s_1 + 1 \geq n$. We distinguish two cases depending on the parity of $|V_2|$.
\medskip

\noindent
Case 1: $|V_2|$ is even.\\
Observe that in this case $a$ is odd. If $2s_1+1 < a$, then $a \geq 2s_1+3$. We can then take the cycle $(s_1+2,s_1+3,\ldots , 2s_1+2)$ of odd length $s_1+1$, the edge $\{1, s_1+1\}$, and cover all of the remaining vertices of $\G'$ by independent $1$-edges, contradicting Assumption~\ref{ass:circ2}. Similarly, if $a - 1 + s_1 > n$, then $j = a - 1 + s_1 - n$ is odd with $1 \leq j \leq s_1-3$ (note that $a+2 \leq n$). We can thus take the cycle $(s_1+1,s_1+2,\ldots , a-1, j, j-1, \ldots , 2,1)$ of odd length $a-s_1+j-1$. This leaves us with the ``last'' $s_1-j-1$ (an even number) vertices of $V_1$ and the whole $V_3$, a contradiction. We therefore find that $2s_1+1 \geq a$ and $a - 1 + s_1 \leq n-1$. Together with the inequality $2a - s_1 + 1 \geq n$, the first of these yields $n \leq 3s_1+3$, while the second yields $n \geq 3s_1-1$. Therefore, $n \in \{3s_1-1, 3s_1+1, 3s_1+3\}$. Moreover, if $n = 3s_1-1$ then $a = 2s_1-1$, while if $n \in \{3s_1+1, 3s_1+3\}$ then $a = 2s_1+1$. 
\medskip

\noindent
Case 2: $|V_2|$ is odd.\\
In this case $a$ is even and $V_3$ (being of odd size) is nonempty. If $2s_1 + 1 < a$, we can take the cycle $(s_1+1, s_1+2, \ldots, 2s_1+1)$ of odd length $s_1+1$, the edge $\{n-1,s_1-1\}$, and cover the remaining vertices of $\G'$ with independent $1$-edges, a contradiction. Similarly, if $a-1+s_1 > n$, then $j = a - 1 + s_1 - n$ is even with $2 \leq j \leq s_1-2$. We can then take the cycle $(1, s_1 + 1, s_1 + 2, \ldots , a - 1, j, j - 1, \ldots , 2)$, the edge $\{s_1-1, n-1\}$, and cover the remaining vertices of $\G'$ by independent $1$-edges, a contradiction. We therefore find that $2s_1 \geq a$ and $a - 1 + s_1 \leq n$. Together with $2a - s_1 + 1 \geq n$, these two inequalities yield $3s_1-3 \leq n \leq 3s_1+1$, implying that $n \in \{3s_1-3, 3s_1-1, 3s_1+1\}$. Moreover, if $n = 3s_1-3$ then $a = 2s_1-2$, while if $n \in \{3s_1-1, 3s_1+1\}$ then $a = 2s_1$. 
\medskip

Recall that we assumed that $|V_2| \geq |V_3|$, with the idea that if this is not the case we can replace the roles of the vertices $0$ and $s_1$. One can verify that there are only two situations, which are not symmetric in the sense that exchanging the roles of the vertices $0$ and $s_1$ changes the corresponding $a$. The first is when $n = 3s_1-1$ and $a = 2s_1$ (where the exchange yields $a = 2s_1-2$), and the other is when $n = 3s_1+1$ and $a = 2s_1+1$ (where the exchange yields $a = 2s_1-1$).
\end{proof}

We finally determine the connection set $S$ and the relationship between $n$ and $s_1$. 

\begin{lemma}
\label{le:circ_x_big_a_even_2}
With reference to Assumption~\ref{ass:circ2}, $n \in \{3s_1-3, 3s_1-1, 3s_1+1, 3s_1+3\}$ and one of the following holds:
\begin{itemize}
\itemsep = 0pt
\item $S = \{\pm 1, \pm s_1\}$,
\item $n = 3s_1-3$ with $s_1 \geq 4$ and $S = \{\pm 1, \pm (s_1-2), \pm s_1\}$,
\item $n = 3s_1+3$ and $S = \{\pm 1, \pm s_1, \pm (s_1+2)\}$.
\end{itemize}
Moreover, in each of these three cases the graph $\G$ is not fractional $2$-extendable.
\end{lemma}

\begin{proof}
By Lemmas~\ref{le:circ_ell_small} and~\ref{le:circ_x_big_a_odd}, $s_1$ is even and $a > s_1$. Moreover,  Lemma~\ref{le:circ_x_big_a_even} implies that $n \in \{3s_1-3, 3s_1-1, 3s_1+1, 3s_1+3\}$ and that one of the items from that lemma holds. Suppose first that $n \in \{3s_1 - 1, 3s_1 + 1\}$ and note that in this case $3$ is coprime to $n$. Multiplication by $3$ is then an isomorphism from $\G$ to the circulant $\Circ(n ; 3S)$. Note that the images of the edges $e$ and $e'$ are the edges corresponding to the elements $1$ and $3$ from $3S$. If $|S| > 4$ we therefore contradict Lemma~\ref{le:circ_x_big_a_odd}. It thus follows that $|S| = 4$, that is, $S = \{\pm 1, \pm s_1\}$. Recall that we already know that $\Circ(n ; \{\pm 1, \pm 3\})$ is not fractional $2$-extendable, implying that $\G$ is also not fractional-$2$-extendable in this case.

From now on we can thus assume that either $n = 3s_1 - 3$ and $a = 2s_1 - 2$, or $n = 3s_1 + 3$ and $a = 2s_1 + 1$. Note that if $s_1 = 2$, then $n = 3s_1+3 = 9$ must hold, and consequently $a = 2s_1+1 = 5$. In this case, $3 \notin S$ since otherwise $(1,3,4,7,8)$ is a spanning $5$-cycle of $\G'$, a contradiction. By Assumption~\ref{ass:circ2} it thus follows that $S = \{\pm 1, \pm 2,\pm 4\}$. By Lemma~\ref{le:CayAb_proper_subgroup}, $\G$ is indeed not fractional $2$-extendable in this case. We can thus assume that $s_1 \geq 4$. Let $V_1$, $V_2$ and $V_3$ be as in the proof of Lemma~\ref{le:circ_x_big_a_odd}. Note that this time $V_2$ and $V_3$ are of equal size $s_1-3$ or $s_1$, depending on whether $n = 3s_1-3$ or $n = 3s_1+3$, respectively. We consider the two possibilities separately.
\medskip

\noindent
Case 1: $n = 3s_1 - 3$.\\
As already mentioned, $a = 2s_1 - 2$ and $|V_2| = |V_3| = s_1 - 3$, which is an odd number. For $s_1 = 4$, we get $n = 9$, in which case multiplication by $2$ gives an isomorphic circulant where the edges $e$ and $e'$ are mapped to edges corresponding to $1$ and $2$ from the new connection set. Since we already settled this possibility, we can thus assume that $s_1 \geq 6$. 
\smallskip

\noindent
We claim that there exists no even $s \in S \setminus \{s_1, s_1-2\}$ with $s < n/2$. Suppose on the contrary that such an $s$ exists. If $s \leq s_1 - 4$, we can take the cycles of (odd) length $s+1$ consisting of the ``first'' $s+1$ vertices in each of $V_1$, $V_2$ and $V_3$, and cover the remaining vertices of $\G'$ by independent $1$-edges, a contradiction. If however $s \geq s_1 + 2$, then $2 \leq 2s_1 - s \leq s_1 - 2$, and we can take the cycle $(n-1,s_1-1,s_1-2,\ldots , 2s_1-s, 2s_1, 2s_1+1, \ldots, n-2)$, the edge $\{1,s_1+1\}$ and cover the remaining vertices of $\G'$ by independent $1$-edges, a contradiction. This proves our claim.
\smallskip

\noindent
We now show there is also no odd $s \in S$ with $1 < s < n/2$. If such an $s$ exists, then $3 \leq s \leq (3s_1-4)/2 = 2s_1-(s_1+4)/2 \leq 2s_1-5$ (recall that $s_1 \geq 6$). But then there exists an odd $j \in V_2$ such that $j+s \in V_3$. We can therefore take the cycle $(1,s_1+1,s_1+2,\ldots , j, j+s, j+s+1, \ldots , n-1, s_1-1, s_1-2, \ldots , 2)$, which leaves us with $2s_1-3-j$ consecutive vertices of $V_2$ and $j+s-2s_1$ consecutive vertices of $V_3$, which can all be covered by independent $1$-edges, a contradiction. 
\smallskip

\noindent
This shows that $S \subseteq \{\pm 1, \pm (s_1-2), \pm s_1\}$. It is easy to see that in this case the graph obtained from $\G'$ by removing its potential edge $\{1, s_1-1\}$ (but leaving its endvertices) is bipartite with $1$ and $s_1-1$ belonging to the smaller set of the bipartition. Therefore, Corollary~\ref{cor:bip} implies that $\G'$ indeed has no fractional perfect matching.
\medskip

\noindent
Case 2: $n = 3s_1+3$.\\
In this case $a = 2s_1+1$ and $|V_2| = |V_3| = s_1$ is even. The argument is very similar to the one in Case 1, so we leave some details to the reader. 
\smallskip

\noindent
We first verify that there exists no even $s \in S \setminus \{s_1, s_1+2\}$ with $s < n/2$. For, if such an $s$ exists, then if $s < s_1$, we can take the cycle containing the ``first'' $s+1$ vertices of $V_1$ and cover the remaining vertices of $\G'$ with independent $1$-edges. If however, $s_1+4 \leq s < n/2$, then $1 \leq 2s_1+3-s \leq s_1-1$, and we can take the cycle $(n-1,s_1-1,s_1-2,\ldots , 2s_1+3-s, 2s_1+3,2s_1+4,\ldots, n-2)$ and cover the remaining vertices of $\G'$ with independent $1$-edges. In both cases we contradict Assumption~\ref{ass:circ2}.
\smallskip

\noindent
Similarly, there exists no odd $s \in S$ with $1 < s < n/2$. For, if it exists, then $3 \leq s \leq 2s_1-1$ (recall that $s_1 \geq 4$), and so there exists an even $j \in V_2$ with $j+s \in V_3$. We then take the cycle $(1,s_1+1,s_1+2, \ldots , j, j+s, j+s+1, \ldots , n-1, s_1-1, s_1-2, \ldots , 2)$ and cover the remaining vertices of $\G'$ with independent $1$-edges, again contradicting Assumption~\ref{ass:circ2}.
\smallskip

\noindent
Therefore, $S \subseteq \{\pm 1, \pm s_1, \pm (s_1+2)\}$. It is again easy to see that the graph obtained from $\G'$ by removing its potential edge $\{s_1+1, n-1\}$ is bipartite with $s_1+1$ and $n-1$ belonging to the smaller set of the bipartition. Therefore, $\G'$ has no fractional perfect matching by Corollary~\ref{cor:bip}.
\end{proof}

Combining together Proposition~\ref{pro:2-f-ext-CayAb_even} and all of the above lemmas, we finally obtain the following classification of fractional $2$-extendable connected Cayley graphs of Abelian groups.

\begin{theorem}
\label{the:CayAb2}
Let $\G = \Cay(A ; S)$ be a connected Cayley graph of an Abelian group of order $n \geq 5$. Then $\G$ is fractional $2$-extendable if and only if it is not isomorphic to one of the following graphs:
\begin{itemize}
\itemsep = 0pt
\item[{\rm (i)}] $\Circ(n;\{\pm 1\})$,
\item[{\rm (ii)}] $\Circ(n; \{\pm 1, 2m\})$, where $n = 4m \geq 8$,
\item[{\rm (iii)}] $\Circ(n; \{\pm 2, 2m+1\})$, where $n = 4m+2 \geq 6$,
\item[{\rm (iv)}] $\Circ(n; \{\pm 1, \pm 2\})$,
\item[{\rm (v)}] $\Circ(n; \{\pm 1, \pm 3\})$, where $n$ is odd,
\item[{\rm (vi)}] $\Circ(n; \{\pm 1, \pm 2m\})$, where $n = 4m+2 \geq 6$,
\item[{\rm (vii)}] $\Circ(n; \{\pm 1, \pm (m-1)\})$, where $n = 3m \geq 9$ with $m$ odd.
\item[{\rm (viii)}] $\Circ(n; \{\pm 1, \pm (m+1)\})$, where $n = 3m \geq 9$ with $m$ odd.
\item[{\rm (ix)}] $\Circ(n; \{\pm 1, \pm (m-1), \pm (m+1)\})$, where $n = 3m \geq 9$ with $m$ odd.
\item[{\rm (x)}] $\Cay(\ZZ_{m} \times \ZZ_3 ; \{\pm (1,0), \pm (1,1)\})$, where $n = 3m \geq 9$ with $m$ odd.
\end{itemize}
\end{theorem}

\section{Concluding remarks}
\label{sec:conclude}

In Subsection~\ref{subsec:even}, we saw that when considering only graphs of even order, the families of fractional and classical $2$-extendable connected Cayley graphs of Abelian groups coincide. Therefore, relaxing the condition of $2$-extendability to that of fractional $2$-extendability does not reduce the class of ``non-examples''. The following question thus arises naturally.

\begin{question}
\label{que1}
Does there exist a Cayley graph of an Abelian group of even order such that for some integer $t$, this graph is fractional $t$-extendable but is not $t$-extendable? If so, what is the smallest $t$ for which such a graph exists?
\end{question}

There are other interesting questions concerning (fractional) extendability in Cayley graphs of Abelian groups. Let us mention just one more. In 1993, Yu~\cite{Yu93} introduced a generalization of the concept of $t$-extendability (in the classical sense) for graphs of even order to those of odd order in the following way. A graph $\G$ of odd order at least $2t+3$ is said to be $t\frac{1}{2}$-extendable (called $t$-near-extendable in~\cite{CioLi14}), if for each vertex $v$ of $\G$ the graph $\G - v$ is $t$-extendable. Using Proposition~\ref{prop:iso}, one can show that a graph which is not fractional $t$-extendable, is also not $t\frac{1}{2}$-extendable, and so each $t\frac{1}{2}$-extendable graph is fractional $t$-extendable. It is thus interesting to investigate the difference between these two concepts, in general, but also in the context of Cayley graphs of Abelian groups of odd order. In~\cite{MikSpa09}, it was shown that, with the exception of cycles of odd length, all connected Cayley graphs of Abelian groups of odd order are $1\frac{1}{2}$-extendable. Corollary~\ref{cor:CayAb1fext} thus implies that within the family of connected Cayley graphs of Abelian groups of odd order, the class of fractional $1$-extendable examples coincides with the class of $1\frac{1}{2}$-extendable ones.  The first next step would be to see if all connected Cayley graphs of Abelian groups of odd order, other than the ones from the ten families from Theorem~\ref{the:CayAb2}, are $2\frac{1}{2}$-extendable. More generally, the following question arises naturally.

\begin{question}
\label{que2}
Does there exist a positive integer $t \geq 2$ and a connected Cayley graph $\G$ of an Abelian group of odd order, such that $\G$ is fractional $t$-extendable, but is not $t\frac{1}{2}$-extendable? If so, what is the smallest $t$ for which such a graph exists?
\end{question}

Finally, it would be interesting to study fractional extendability (and the corresponding analogues of Questions~\ref{que1} and~\ref{que2}) for some other nice families of graphs. One could consider Cayley graphs of other groups or vertex-transitive graphs in general. However, other families of graphs with a large degree of regularity could also be of interest. Examples of such families, where investigations of classical extendability has been considered, are edge-regular graphs~\cite{KutMarMikSpa19}, quasi-strongly regular graphs~\cite{AlaHusMikSpa18}, strongly regular graphs~\cite{CioLi14, HL, LZ96}, and more generally distance regular graphs~\cite{CioKooLi17}, to mention just a few. The topic thus offers many possible directions for future investigations.

\section*{Acknowledgments}

The authors would like to thank the anonymous referees for helpful suggestions to improve the paper, and to \v Stefko Miklavi\v c for introducing them to the concept of fractional extendability and fruitful discussions on the problems considered. They also acknowledge the financial support by the Slovenian Research and Innovation Agency (research program P1-0285 and research projects J1-3001 and J1-50000).

\end{document}